\documentclass[11pt]{article}
\usepackage{amsmath,amsfonts,geompsfi}

\renewcommand{\bold}[1]{\medskip \noindent {\bf \boldmath #1
                        }\nopagebreak[4]}

\newtheorem{theorem}{Theorem}[section]
\newcommand{\cx}{{\mathbb C}}
\newcommand{\half}{{\mathbb H}}
\newcommand{\integers}{{\mathbb Z}}

\newcommand{\reals}{{\mathbb R}}

\newcommand{\makefig}[3]{
	\begin{figure}[htbp]
        \refstepcounter{figure}
	\label{#2}
        \begin{center}~
		#3~\\
		\medskip
                {\sf Figure \thefigure.  #1}
        \end{center}
	\medskip
	\end{figure}
}

\newenvironment{pf*}[1]{%
 \begin{proof}[#1]%
}{ 
 \end{proof}
}

\newcommand{\ital}[1]{\medskip \noindent {\em #1 }\nopagebreak[4]}

\newcommand{\bdry}{\partial}

\newcommand{\closure}{\overline}

\newcommand{\del}{\partial}

\newcommand{\disjunion}{\sqcup}

\newcommand{\includesin}{\hookrightarrow}

\usepackage{amssymb}

\newcommand{\st}{\; | \;}         %Such that

\newcommand{\wt}{\widetilde}

\newcommand{\zed}{\integers}
\newcommand{\inj}{\mbox{\rm inj}}

\newcommand{\interior}{\mbox{\rm int}}

\newcommand{\Isom}{\mbox{\rm Isom}}

\newcommand{\PSL}{\mbox{\rm PSL}}

\newcommand{\Teich}{\mbox{\rm Teich}}

\newtheorem{conj}[theorem]{Conjecture}
\newtheorem{defn}[theorem]{Definition}

\newcommand{\calA}{{\mathcal A}}

\newcommand{\calL}{{\mathcal L}}
\newcommand{\calM}{{\mathcal M}}

\newcommand{\calP}{{\mathcal P}}

\newcommand{\pl}{{\calP \calL}}

\usepackage{euscript}

\newcommand{\eM}{{\EuScript M}}

\newcommand{\bld}{\bold}
\title{\vspace{-.5in}
{\bf Cone-manifolds and the\\ density conjecture}
}
\author{Jeffrey F. Brock\thanks{Research supported by NSF
research grants and an NSF postdoctoral fellowship.}\ \ and Kenneth
W. Bromberg\thanks{Research 
supported by NSF research grants and the Clay Mathematics
Institute. \hfill
}}
\date{June 2, 2002}

\begin{document}

\maketitle
\begin{abstract}
We give an expository account of our proof that each cusp-free
hyperbolic 3-manifold $M$ with finitely generated fundamental group and
incompressible ends is an algebraic limit of geometrically finite
hyperbolic 3-manifolds.
\end{abstract}

\section{Introduction}

The aim of this paper is to outline and describe new constructions and
techniques we hope will provide a useful tool to study
deformations of hyperbolic 3-manifolds.  An initial application
addresses the following conjecture.
\begin{conj}[Bers-Sullivan-Thurston]{\sc The Density Conjecture}
Each complete hyperbolic 3-manifold $M$ with finitely generated
fundamental group is an algebraic limit of geometrically finite
hyperbolic 3-manifolds.
\label{conjecture:density}
\end{conj}
Algebraic convergence of $M_n$ to $M$ refers to convergence in the
{\em algebraic deformation space} or in the topology of convergence on
generators of the holonomy representations $$\rho_n \colon \pi_1(M)
\to \PSL_2(\cx) =
\Isom^+ (\half^3).$$  The approximating manifolds $M_n = \half^3
/\rho_n(\pi_1(S))$ are {\em geometrically finite} if the {\em convex
core} of $M_n$, the minimal
convex subset homotopy equivalent to $M_n$, has finite volume.
We give an expository account of our progress toward
Conjecture~\ref{conjecture:density} \cite{Brock:Bromberg:density}.
\begin{theorem}
Let $M$ be a complete hyperbolic 3-manifold with no cusps, finitely
generated fundamental group, and incompressible ends.  Then $M$ is an
algebraic limit of geometrically finite hyperbolic 3-manifolds.
\label{theorem:main}
\end{theorem}
Our result represents an initial step in what we hope will be a
general geometrically finite approximation theorem for {\em
topologically tame} complete hyperbolic 3-manifolds, namely, for each
such manifold $M$ that is homeomorphic to the interior of a compact
3-manifold.  
%Such a result would reduce
%Conjecture~\ref{conjecture:density} to Marden's original conjecture.
%\begin{conj}[Marden]
%\label{conjecture:tame}
%Each complete hyperbolic 3-manifold with finitely generated
%fundamental group is topologically tame.
%\end{conj}

Indeed, the clearly essential assumption in our argument is that 
$M$ is tame; we make direct use of the following theorem due to Bonahon and
Thurston.
\begin{theorem}[Bonahon-Thurston]
Each cusp-free complete hyperbolic 3-manifold $M$ with finitely
generated fundamental group and incompressible
ends is geometrically and topologically tame.
\label{theorem:tameness}
\end{theorem}
(See \cite{Bonahon:tame, Thurston:book:GTTM}).  The tameness of a
complete hyperbolic 3-manifold with finitely generated fundamental
group reduces to a consideration of its {\em ends} since every such
3-manifold $M$ contains a {\em compact core}, namely, a compact
submanifold $\calM$ whose inclusion is a homotopy equivalence.  Each
end $e$ of $M$ is associated to a component $E$ of $M \setminus
\interior(\calM)$, which we typically refer to as an ``end'' of $M$,
assuming an implicit choice of compact core.  An end $E$ is {\em
incompressible} if the inclusion of $E$ induces an injection $\pi_1(E)
\includesin \pi_1(M)$.  The end $E$ is {\em geometrically finite} if
it has compact intersection with the convex core.  Otherwise, it is
{\em degenerate}.

For a degenerate end $E$, {\em geometric} tameness refers to the
existence of a family of simple closed curves on the closed surface $S
= \bdry \calM \cap 
E$ whose geodesic representatives leave every compact subset of $E$.
Using interpolations of pleated surfaces, Thurston showed that a
geometrically tame end is homeomorphic to $S \times \reals^+$, so
$M$ is topologically tame if all its ends are geometrically finite or
geometrically tame (R. Canary later proved the equivalence of these
notions \cite{Canary:ends}).

\bld{Approximating the ends.}  Our approach to
Theorem~\ref{theorem:main} will be to approximate the manifold $M$
{\em end by end}.  Such an approach is justified by an {\em asymptotic
isolation} theorem (Theorem~\ref{theorem:isolation:intro}) that
isolates the geometry of the ends of $M$ from one another when $M$ is
obtained as a limit of geometrically finite manifolds. 
Each degenerate end $E$ of $M$ has one of two types: $E$ has either
\begin{enumerate}
\item[{\bf I.}] {\sl bounded geometry:} there is a uniform lower bound to the
length of the shortest geodesic in $E$, or
\item[{\bf II.}] {\sl arbitrarily short geodesics:} there is some sequence
$\gamma_n$ of geodesics in $E$ whose length is tending to zero.
\end{enumerate}
Historically, it is the latter category of ends that have been persistently
inscrutable (they are known to be generic
\cite{McMullen:cusps,CCHS:density}).  
Our investigation of such ends
begins with another key consequence of tameness, due to J. P. Otal
(see \cite{Otal:unknotted}, or his article \cite{Otal:unknotted2} in
this volume). 
Before discussing this result, we introduce some terminology.

If $E$ is an incompressible end of $M$, the cover $\widetilde{M}$
corresponding to $\pi_1(E)$ is homotopy equivalent to the surface $S =
\bdry \calM \cap E$.  Thus, $\wt M$ sits in the {\em algebraic
deformation space} $AH(S)$, namely, hyperbolic 3-manifolds $M$
equipped with homotopy equivalences, or {\em markings}, $f \colon S
\to M$ up to isometries that preserve marking and orientation  (see
\cite{Thurston:hype2}, 
\cite{McMullen:book:RTM}).  The space $AH(S)$ is equipped with the
algebraic topology, or the topology of convergence of holonomy
representations, as described above.  Theorem~\ref{theorem:tameness}
guarantees each $M \in AH(S)$ is homeomorphic to $S \times \reals$;
Otal's theorem provides deeper information about how short geodesics
in $M$ sit in this product structure.
\begin{theorem}[Otal \cite{Otal:unknotted}]
Let $M$ lie in $AH(S)$.  There is an $\epsilon_{\rm knot} >0$ so that if
$\calA$ is any collection of closed geodesics so that for each $\gamma
\in \calA$ we have $$\ell_M(\gamma) < \epsilon_{\rm knot}$$ then there exists a
collection of distinct real numbers $\{t_\gamma \st \gamma \in
\calA\}$ and an ambient isotopy of $M \cong S \times \reals$ taking
each $\gamma$ to a simple curve in $S \times \{t_\gamma\}$.
\label{theorem:Otal}
\end{theorem}
Said another way, sufficiently short curves in $M$ are simple,
unknotted and pairwise unlinked with respect to the product structure
$S \times \reals$ on $M$.

Otal's theorem directly facilitates the {\em grafting} of tame ends
that carry sufficiently short geodesics.  This procedure, introduced
in \cite{Bromberg:bers}, uses embedded end-homotopic annuli in a
degenerate end to perform 3-dimensional version of grafting from the
theory of projective structures (see e.g. \cite{McMullen:graft,Gallo:Kapovich:Marden}).  In section~\ref{section:grafting} we
will describe how successive graftings about short curves in an end $E$
of $M$ can be used to produce a sequence of projective structures with
holonomy $\pi_1(M)$ whose underlying conformal structures $X_n$
reproduce the asymptotic geometry of the end $E$ in a limit.

Our discussion of ends $E$ with bounded geometry relies directly on a
large body of work of Y. Minsky
\cite{Minsky:Teichmuller,Minsky:ends,Minsky:KGCC,Minsky:bounded} which
has recently resulted in the following {\em bounded geometry theorem}.
\begin{theorem}[Minsky]{\sc Bounded Geometry Theorem}
Let $M$ lie in $AH(S)$, and assume $M$ has a global lower bound
to its injectivity radius $\inj \colon M \to \reals^+$.  If $N \in
AH(S)$ has the same end-invariant as that of $M$ then $M=N$ in
$AH(S)$.
\label{theorem:bounded}
\end{theorem}
In other words, there is an orientation preserving isometry $\varphi
\colon M \to N$ that respects the homotopy classes of the markings on
each.  The ``end invariant'' $\nu(M)$ refers to a union of invariants,
each associated to an end $E$ of $M$.  Each invariant is either a
Riemann surface in the {\em conformal boundary} $\bdry M$ that
compactifies the end, or an {\em ending lamination}, namely, the
support $|\mu|$ of a limit $[\mu]$ of simple closed curves $\gamma_n$
whose geodesic representatives in $M$ that exit the end $E$ (here
$[\mu]$ is the limit of $[\gamma_n]$ in Thurston's projective measured
lamination space $\pl(S)$ \cite{Thurston:book:GTTM, Thurston:hype2}).

Minsky's theorem proves Theorem~\ref{theorem:main} for each $M$ with a
lower bound to its injectivity radius, since given any end invariant
$\nu(M)$ there is some limit $M_\infty$ of geometrically finite manifolds with
end invariant $\nu(M_\infty) = \nu(M)$ (see \cite{Ohshika:ends,
Brock:length}). 

\bld{Realizing ends on a Bers boundary.}  Grafting ends with short
geodesics and applying Minsky's results to ends with bounded geometry,
we arrive at a {\em realization} theorem for ends of manifolds $M \in
AH(S)$ in some {\em Bers compactification}.
\begin{theorem}{\sc Ends are Realizable}
Let $M \in AH(S)$ have no cusps.  Then each end of $M$ is
realized in a Bers compactification.
\label{theorem:ends:realizable}
\end{theorem}
We briefly explain the idea and import of the theorem.  The
subset of $AH(S)$ consisting of geometrically finite cusp-free
manifolds is the {\em quasi-Fuchsian locus} $QF(S)$.  In
\cite{Bers:simunif} Bers exhibited
the parameterization $$Q \colon \Teich(S)
\times \Teich(S) \to QF(S)$$
so that $Q(X,Y)$ contains $X$ and $Y$ in its conformal boundary;
$Q(X,Y)$ {\em simultaneously uniformizes} the pair $(X,Y)$.  Fixing
one factor, we obtain the {\em
Bers slice} $B_Y = \{ Q(X,Y) \st Y \in \Teich(S) \}$, which Bers proved
to be precompact.  The resulting compactification $\closure{B_Y}
\subset AH(S)$ for
Teichm\"uller space has frontier $\bdry B_Y$, a {\em Bers boundary}
(see \cite{Bers:bdry}).

We say an end $E$ of $M \in AH(S)$ is {\em realized by $Q$ in the
Bers compactification $\closure{B_Y}$} if there is a manifold $Q
\in \closure{B_Y}$ 
and a marking preserving bi-Lipschitz embedding $\phi \colon E \to Q$
(see Definition~\ref{definition:ends:realizable}).  

The cusp-free manifold $M \in AH(S)$ is {\em singly-degenerate} if
exactly one end of $M$ is compactified by a conformal boundary
component $Y$.  In this case, the main theorem of \cite{Bromberg:bers}
establishes that $M$ itself lies in the Bers boundary $\bdry B_Y$,
which was originally conjectured by Bers \cite{Bers:bdry}.
Theorem~\ref{theorem:ends:realizable} generalizes this result to the
relative setting of a given incompressible end of $M$, allowing us to
pick candidate approximates for a given $M$ working end-by-end.

\bld{Candidate approximates.} To see explicitly how candidate
approximates are chosen, let $M$ have finitely generated fundamental
group and incompressible ends.  For each end $E$ of $M$,
Theorem~\ref{theorem:ends:realizable} allows us to choose $X_n(E)$ so
that the limit of $Q(X_n(E),Y)$ in $\closure{B_Y}$ realizes the end
$E$.  Then we simply let $M_n$ be the geometrically finite 
manifold homeomorphic to $M$ determined by specifying the data $$(X_n(E_1),
\ldots, X_n(E_m)) \in \Teich(\bdry \calM)$$ where $\calM$ is a compact
core for $M$; $\Teich(\bdry \calM)$ naturally parameterizes such
manifolds (see section~\ref{section:general}).  The union $X_n(E_1)
\cup \ldots \cup X_n(E_m)$ constitutes the conformal boundary $\bdry
M_n$.

To conclude that the limit of $M_n$ is the original manifold
$M$, we must show that limiting geometry of each end of $M_n$ does not
depend on limiting phenomena in the other ends.  
We show ends of $M_n$ are {\em asymptotically isolated}.
\begin{theorem}{\sc Asymptotic Isolation of Ends} 
\label{theorem:isolation:intro}
Let $N$ be a complete cusp-free hyperbolic 3-manifold with  finitely generated
fundamental group and incompressible ends.
Let $M_n$ converge algebraically to $N$.  Then
up to bi-Lipschitz diffeomorphism, the end $E$ of $M$
depends only on the corresponding sequence $X_n(E) \subset \bdry M_n$.
\end{theorem}
(See Theorems~\ref{theorem:isolation} and~\ref{theorem:isolation:II} 
for a more precise formulation).

When $N \in AH(S)$ is singly-degenerate, the theorem is well known
(see, e.g. \cite[Prop. 3.1]{McMullen:book:RTM}).
For $N$ not homotopy equivalent to a surface, the cover corresponding
to each end of $N$ is singly-degenerate, so the theorem follows in
this case as well.

The ideas in the proof of Theorem~\ref{theorem:isolation:intro} when
$N$ is doubly-degenerate represent a central focus of this paper.  In
this case, the cover of $N$ associated to each end is again the
manifold $N$ and thus not singly-degenerate, so the asymptotic
isolation is no longer immediate.  The situation is remedied by a new
technique in the cone-deformation theory called the {\em drilling
theorem} (Theorem~\ref{theorem:bilip}).

This drilling theorem allows us to ``drill out'' a sufficiently short
curves in a geometrically finite cusp-free manifold with bounded
change to the metric outside of a tubular neighborhood of the drilling
curve.  When quasi-Fuchsian manifolds $Q(X_n,Y_n)$ converge to the
cusp-free limit $N$, any short geodesic $\gamma$ in $N$ may be drilled
out of each $Q(X_n,Y_n)$.  

The resulting drilled manifolds
$Q_n(\gamma)$ converge to a limit $N(\gamma)$ whose higher genus ends
are bi-Lipschitz diffeomorphic to those of $N$.  In the manifold
$N(\gamma)$, the rank-2 cusp along $\gamma$ serves to insulate the
geometry of the ends from one another, giving the necessary
control. (When there are no short curves, Minsky's theorem again
applies).

The drilling theorem manifests the idea that the thick part of a
hyperbolic 3-manifold with a short geodesic looks very similar to the
thick part of the hyperbolic 3-manifold obtained by removing that
curve.  We employ the cone-deformation theory of C. Hodgson and
S. Kerckhoff to give analytic control to this qualitative picture.

\bld{Plan of the paper.}  In what follows we will give descriptions
of each facet of the argument.  Our descriptions are expository in
nature, in the interest of conveying the main ideas rather than
detailed specific arguments (which appear in
\cite{Brock:Bromberg:density}).  We will focus on the case when $M$ is
homotopy equivalent to a surface, which presents the primary
difficulties, treating the general case briefly at the conclusion.

In section~\ref{section:cone} we provide an overview of techniques in
the deformation theory of hyperbolic cone-manifolds we will apply,
providing bounds on the metric change outside a tubular
neighborhood of the cone-singularity under a change in the cone-angle.
In section~\ref{section:grafting} we describe the grafting
construction and how it produces candidate approximates for the ends
of $M$ with arbitrarily short geodesics.  
Section~\ref{section:drilling} describes the asymptotic isolation
theorem (Theorem~\ref{theorem:isolation:intro}), the realization
theorem for ends (Theorem~\ref{theorem:ends:realizable}),
and finally how these results combine to give a proof of
Theorem~\ref{theorem:main} when $M$ lies in $AH(S)$.
The general case is discussed in section~\ref{section:general}.

\bld{Acknowledgments.}  The authors 
would like to thank Craig Hodgson and Steve Kerckhoff for their
support and for providing much of the analytic basis for our results,
Dick Canary and Yair Minsky for their input and inspiration, and
Caroline Series for her role in organizing the 2001 Warwick conference
and for her solicitation of this article.

\section{Cone-deformations}
\label{section:cone}
Over the last decade, Hodgson and Kerckhoff have developed a powerful
rigidity and deformation theory for 3-dimensional 
hyperbolic cone-manifolds 
\cite{Hodgson:Kerckhoff:rigidity}.  While their theory was developed
initially for application to closed hyperbolic cone-manifolds,
work of the second author (see \cite{Bromberg:thesis}) has generalized
this rigidity and deformation theory to {\em infinite volume}
geometrically finite manifolds.  

The cone-deformation theory represents a key technical tool in
Theorem~\ref{theorem:main}. Let $N$ be a compact, hyperbolizable
3-manifold with boundary; assume that $\del N$ does not contain tori
for simplicity. Let $c$ be a simple closed curve in the interior of
$N$. A {\em hyperbolic cone-metric} is a hyperbolic metric on the
interior of $N\setminus c$ that completes to a singular metric on all of the
interior of $N$. Near $c$ the metric has the form $$dr^2 + \sinh^2 r
d\theta^2 + \cosh^2 r dz^2$$ where $\theta$ is measured modulo the
{\em cone-angle}, $\alpha$.

Just as $\half^3$ is compactified by the Riemann sphere, complete
infinite volume hyperbolic 3-manifolds are often compactified by
projective structures. If a hyperbolic cone-metric is so compactified
it is {\em geometrically finite without rank-one cusps}. As we have
excised the presence of rank-one cusps in our hypotheses, we simply
refer to such metrics as geometrically finite.

A projective structure on $\del N$ has an underlying conformal
structure; we often refer to $\partial N$ together with its conformal
structure as the {\em conformal boundary of $N$}.  
\begin{theorem}
\label{param}
Let $M_\alpha$ denote $N$ with a 3-dimensional geometrically finite
hyperbolic cone-metric with cone-angle $\alpha$ at $c$. If the
cone-singularity has tube-radius at least $\sinh^{-1} 
\left(\sqrt{2}\right)$, then nearby cone-metrics are locally
parameterized by the cone-angle and the conformal boundary.
\end{theorem}
Here, the {\em tube-radius} about  $c$ is the
radius of the maximally embedded metric tube about $c$ in $M_\alpha$.

This local parameterization theorem was first proven by Hodgson and
Kerckhoff for closed manifolds with cone-angle less than $2 \pi$ and
no assumption on the size of the tube radius
\cite{Hodgson:Kerckhoff:rigidity}.  In the thesis of the second author
\cite{Bromberg:thesis}, Hodgson and Kerckhoff's result was generalized
to the setting of general geometrically finite cone-manifolds, where
the conformal boundary may be non-empty.  The replacement of the
cone-angle condition with the tube-radius condition is recent work of
Hodgson and Kerckhoff \cite{Hodgson:Kerckhoff:shape}.

Theorem~\ref{param} allows us to decrease the cone-angle while
keeping the conformal boundary fixed at least for cone-angle near
$\alpha$. We need more information if we wish to decrease the
cone-angle all the way to zero. 
\begin{theorem}[\cite{Bromberg:Schwarzian}]
\label{downtozero}
Let $M_\alpha$ be a 3-dimensional geometrically finite hyperbolic
cone-metric with cone-angle $\alpha$. Suppose that the
cone-singularity $c$ has tube-radius at least $ \sinh^{-1} 
\left(\sqrt{2}\right)$. Then there exists an $\epsilon > 0$ depending
only on $\alpha$ such that if the length of $c$ is less than
$\epsilon$ there exists a one-parameter family, $M_t$, of
geometrically finite cone-metrics with cone-angle $t$ and conformal
boundary fixed for all $t \in [0, \alpha]$.
\end{theorem}

\bld{The drilling theorem.}  When the cone-angle $\alpha$ is $2\pi$ the
hyperbolic cone-metric $M_\alpha$ is actually a smooth hyperbolic
metric. When the cone-angle is zero the hyperbolic cone-metric is also
a smooth complete metric; the curve $c$, however, has receded to
infinity leaving a rank-two cusp, and the complete hyperbolic metric
lives on the interior of $N \setminus c$.  We call $N \setminus c$
with its complete hyperbolic metric $M_0$ the {\em drilling along $c$} of
$M_\alpha$.

Applying analytic tools and estimates developed by Hodgson and
Kerckhoff \cite{Hodgson:Kerckhoff:bounds}, we obtain infinitesimal
control on the metric change outside a tubular neighborhood of the
cone-singularity under a change in the cone-angle.  Letting $U_t
\subset M_t$ denote a standard tubular neighborhood of the
cone-singularity  
we obtain the following {\em
drilling theorem}, which summarizes the key geometric information
emerging from these estimates.
\begin{theorem}{\sc The Drilling Theorem}
\label{theorem:bilip}
Suppose $M_\alpha$ is a geometrically finite hyperbolic cone-metric
satisfying the conditions of Theorem~\ref{downtozero}, and let $M_t$ be
the resulting family of cone-metrics. Then for each $K > 1$ there
exists an $\epsilon'>0$ depending only on $\alpha$ and $K$ such that if
the length of $c$ is less than $\epsilon'$, there are diffeomorphisms
of pairs
$$\phi_t:(M_\alpha \setminus U_\alpha, \bdry U_\alpha)
\longrightarrow (M_t \setminus U_t,\bdry U_t)$$
so that $\phi_t$ is $K$-bi-Lipschitz for each $t \in [0,\alpha]$,
and $\phi_t$ extends over $U_\alpha$ to a homeomorphism for each 
$t \in (0,\alpha]$. 
\end{theorem}

\section{Grafting short geodesics}
\label{section:grafting}

A simple closed curve $\gamma$ in $M \in AH(S)$ is {\em unknotted} if
it is isotopic in $M$ to a simple curve $\gamma_0$ in the ``level
surface'' $S
\times \{ 0 \}$ in the product structure $S \times \reals$ on $M$.
For such a $\gamma$, there is a bi-infinite
annulus $A$ containing $\gamma$ representing its free homotopy class
so that $A$ is
isotopic to $\gamma_0 \times \reals$.  Let $A^+$ denote the sub-annulus of
$A$  exiting the positive end of $M$, let $A^-$ denote the sub-annulus
of $A$ exiting the negative end.

The {\em positive grafting ${\rm Gr}^+(\gamma,M)$ of $M$ along $\gamma$} is
the following surgery of $M$ along the {\em positive grafting annulus} $A^+$.  
\begin{enumerate}

\item Let $M_\zed$ denote the cyclic cover of $M$ associated to the
curve $\gamma$.  Let $$F \colon S^1 \times [0,\infty) \to A^+$$ be a
parameterization of the grafting annulus and let $F_\zed$ be its lift
to $M_\zed$.
\makefig{The grafting annulus and its
lift.}{figure:annulus}{\psfig{file=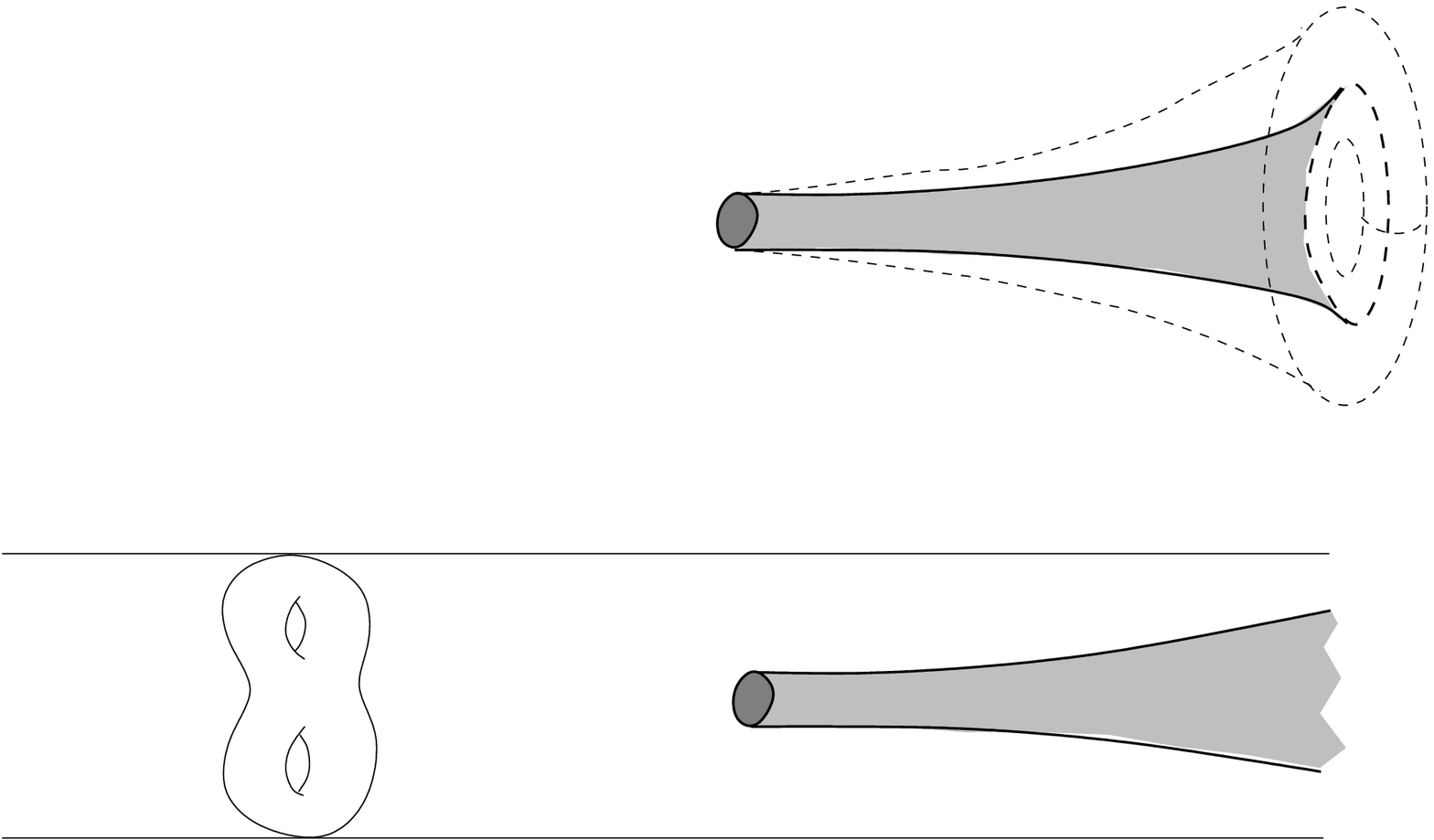,height=2.5in}} 

\item Cutting $M$ along $A^+$ and $M_\zed$ along $A_\zed^+ = F_\zed(S^1
\times [0,1))$,
the complements $M \setminus A^+$ and $M_\zed \setminus
A_\zed^+$ each have two isometric copies of the annulus in their metric
completions $\closure{M \setminus A^+}$ and $\closure{M_\zed
\setminus A_\zed^+}$: the {\em inward annulus} inherits an orientation from $F$
that agrees with the orientation induced by the positive orientation
on $M \setminus A^+$ and the {\em outward annulus} inherits the
opposite orientations from $F$ and $M
\setminus A^+$.  The complement $M_\zed \setminus
A_\zed^+$ also contains an {\em inward} and {\em outward} copy
of $A_\zed^+$ in its metric completion.  

\item Let $F^{\rm in}$ and
$F^{\rm out}$ denote the natural parameterizations of the inward and
outward annulus for the metric completion of $M \setminus A^+$ induced
by $F$ and let $F_\zed^{\rm in}$ and ${F}_\zed^{\rm
out}$ be similarly induced by ${F}_\zed$.

\item Let
$\phi$ be the isometric gluing of 
the inward annulus for $\closure{ M_\zed \setminus A_\zed^+}$ 
to 
the outward annulus for $\closure{M \setminus A^+}$ 
and 
the outward annulus of $\closure{M_\zed \setminus {A_\zed^+}}$ 
to
the inward annulus of $\closure{M \setminus A^+}$ 
so that
$$\phi(F^{\rm in}(x,t)) = {F}_\zed^{\rm out}(x,t) \ \ \ \text{and} \ \
\ \phi(F^{\rm out}(x,t)) = {F}_\zed^{\rm in}(x,t)$$ (the map $\phi$ on
the geodesic $\widetilde{\gamma} \subset M_\zed$ should just be the
restriction covering map $M_\zed \to M$).

\makefig{Grafting: glue the wedge $M_\zed \setminus 
A_\zed^+$ along the completion of $M\setminus
A^+$.}{figure:graft}{\psfig{file=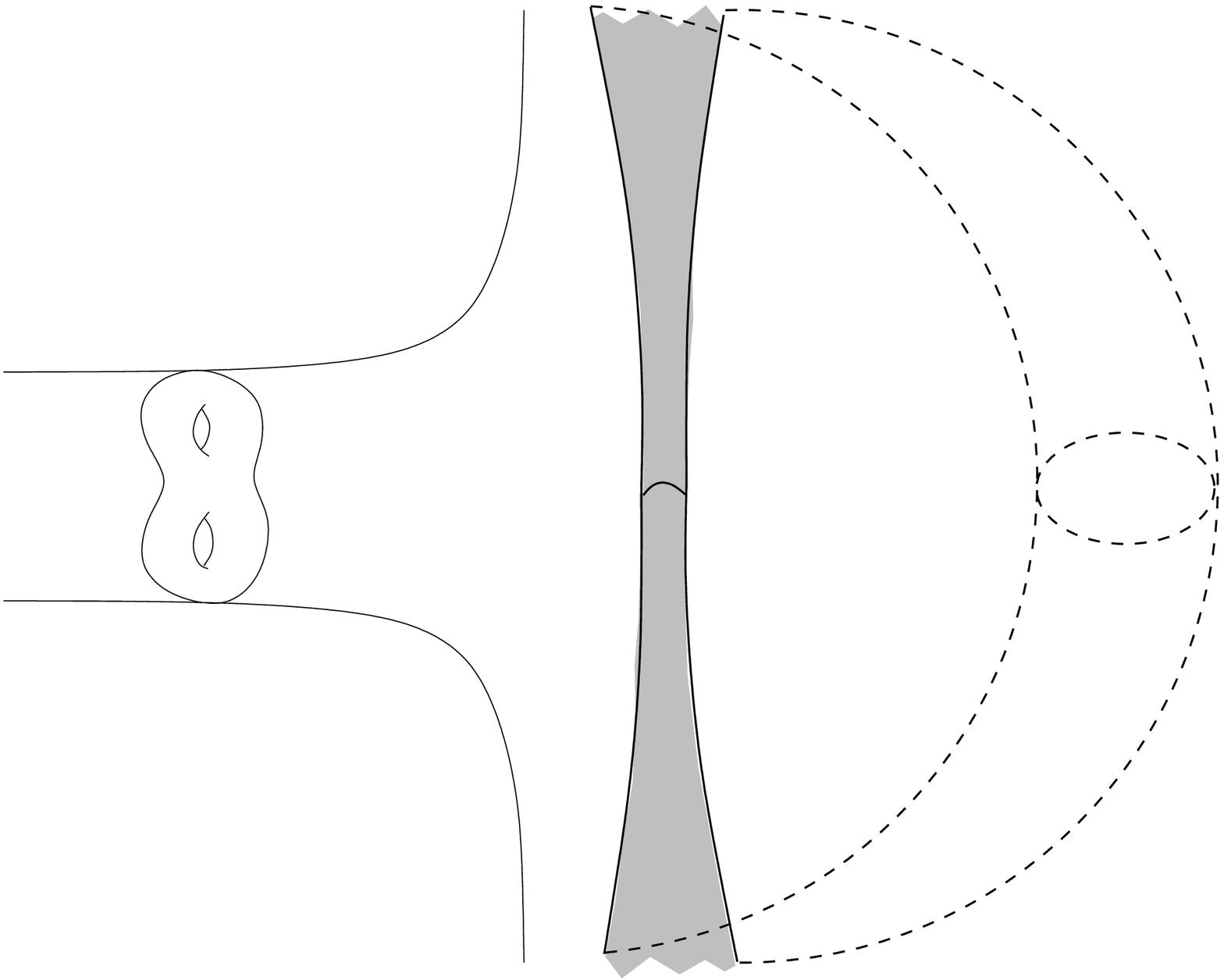,height=2.5in}} 
\end{enumerate}
The result ${\rm Gr}^+(M,\gamma)$ of positive grafting along $\gamma$
is no longer a smooth manifold since its metric is not smooth at
$\gamma$, but ${\rm Gr}^+(M,\gamma)$ inherits a smooth hyperbolic 
metric from $M$ and $M_\zed$ away from $\gamma$.

\bld{Graftings as cone-manifolds.}
Otal's theorem (Theorem~\ref{theorem:Otal}) guarantees that a
sufficiently short closed geodesic $\gamma^*$ is unknotted.
In this case, the positive grafting ${\rm Gr}^+(M,\gamma^*)$ along the
closed geodesic $\gamma^*$ is well defined, and the singularity has a
particularly nice structure: since the singularity is a geodesic, the
smooth hyperbolic structure on ${\rm Gr}^+(M,\gamma^*) \setminus
\gamma^*$ extends to a hyperbolic cone-metric
on ${\rm Gr}^+(M,\gamma^*)$ with cone-singularity $\gamma^*$ and
cone-angle $4 \pi$ at $\gamma^*$ (cf. \cite{Bromberg:bers}).

\bld{Simultaneous grafting.}  We would like to apply the
cone-deformation theory of section~\ref{section:cone} to the grafting
${\rm Gr}^+(M,\gamma^*)$.  The deformation theory applies, however,
only to geometrically finite hyperbolic cone-manifolds.  The grafting
${\rm Gr}^+(M,\gamma^*)$ alone may not be geometrically finite if the
manifold $M$ is doubly-degenerate.  Indeed, in the doubly-degenerate
case positive grafting produces a geometrically finite {\em positive}
end, but to force geometric finiteness of {\em both} ends, we must
perform negative grafting as well.

Let $\gamma$ and $\beta$ be two simple unknotted curves in $M$ that
are also {\em unlinked}: $\gamma$ is isotopic to a level surface in
the complement of $\beta$.  Then $\gamma$ is homotopic either to
$+\infty$ or to $-\infty$ in the complement of $\beta$.  Assume the
former.  Then we may choose a positive grafting annulus $A_\gamma^+$
for $\gamma$ and a negative grafting annulus $A_\beta^-$ for $\beta$
and perform {\em simultaneous grafting} on $M$: we simply perform the
grafting surgery on $A_\gamma^+$ and $A_\beta^-$ at the same time.

By Otal's theorem, when $\gamma^*$ and $\beta^*$ are sufficiently
short geodesics in the hyperbolic 3-manifold $M$, they are simple,
unknotted and unlinked.  If $\gamma^*$ is homotopic to $+\infty$ in
$M \setminus \beta^*$, 
the simultaneous grafting
$${\rm Gr}^\pm(\beta^*, \gamma^* ,M)$$ produces a hyperbolic
cone-manifold with two cone-singularities, one at $\gamma^*$ and one
at $\beta^*$, each with cone-angle $4\pi$.

We then prove the following theorem.
\begin{theorem}{\sc Simultaneous Graftings}  Let $\gamma^*$ and
$\beta^*$ be two simple closed geodesics in $M$ with $\gamma^*$
isotopic to $+\infty$ in the complement of $\beta^*$.  Then the
simultaneous grafting ${\rm Gr}^\pm(\beta^*,\gamma^*,M)$ is a
geometrically finite hyperbolic cone-manifold.
\label{theorem:finite}
\end{theorem}

The proof applies the theory of geometric finiteness for variable
negative curvature developed by Brian Bowditch \cite{Bowditch:hulls}
\cite{Bowditch:finiteness}, to a variable negative curvature smoothing
$\eM$ of ${\rm Gr}^\pm(\beta^*,\gamma^*,M)$ at its cone-singularities.
Using these results, we obtain the following version of Canary's
geometric tameness theorem \cite{Canary:ends} for Riemannian
3-manifolds with curvature pinched between two negative constants, or
{\em pinched negative curvature} (we omit the cusped case as usual).
\begin{theorem}{\sc Geometric Tameness for Negative Curvature}
Each end $E$ of the topologically tame 3-manifold $\eM$ with pinched
negative curvature and no cusps satisfies the following
dichotomy: either
\begin{enumerate}
\item {\sl $E$ is geometrically finite:} $E$ has finite volume
intersection with the convex core of $\eM$, or
\item {\sl $E$ is simply degenerate:} there are essential,
non-peripheral simple closed curves
$\gamma_n$ on the surface $S$ cutting off $E$ whose geodesic
representatives exit every compact subset of $E$.
\end{enumerate}
\end{theorem}

In our setting, any simple closed curve $\eta$ on $S$
whose geodesic representative $\eta^*$ avoids the cone-singularities of
 ${\rm Gr}^\pm(\beta^*, \gamma^*, M)$ 
projects to a closed geodesic 
$\pi(\eta^*)$ in $M$ under the natural local isometric covering $$\pi
\colon {\rm Gr}^\pm(\beta^*, \gamma^*, M) \setminus \beta^* \disjunion
\gamma^* \to M.$$  

The projection $\pi$ extends to a homotopy equivalence across $\beta^*
\disjunion \gamma^*$, so the image $\pi(\eta^*)$ is the geodesic
representative of $\eta$ in $M$.  Though $\pi$ is not proper, we show
that any sequence $\eta_n$ of simple closed curves on $S$ whose
geodesic representatives in ${\rm Gr}^\pm(\beta^*, \gamma^*, M)$ leave
every compact subset must have the property that $\pi(\eta_n^*)$
leaves every compact subset of $M$.  This contradicts bounded diameter
results from Thurston's theory of pleated surfaces
\cite{Thurston:book:GTTM}, which guarantee 
that realizations of $\pi(\eta_n^*)$ by pleated surfaces remain in a
compact subset of $M$.  The contradiction implies that grafted ends
are geometrically finite, proving Theorem~\ref{theorem:finite}.

The simultaneous grafting ${\rm Gr}^+(\beta^*,\gamma^*,M)$ has two
components in its {\em projective boundary at infinity} to which the
hyperbolic cone-metric extends.  Already, we can give an outline of
the proof of Theorem~\ref{theorem:main} in the case that each end of
the doubly-degenerate manifold $M \in AH(S)$ has arbitrarily short
geodesics.

Here are the steps:
\begin{enumerate}
\item[{\bf I.}]  Let $\{\gamma_n^*\}$ be arbitrarily short geodesics
exiting the positive end of $M$ and let $\{\beta_n^*\}$ be arbitrarily
short geodesics exiting the negative end of $M$.  Assume $\gamma_n^*$
is homotopic to $+\infty$ in $M \setminus \beta_n^*$.

\item[{\bf II.}]  The simultaneous graftings 
$${\rm Gr}^\pm(\beta_n^*,\gamma_n^*,M) = M_n^c$$ have projective
boundary with underlying conformal structures $X_n$ on the negative
end of $M_n^c$ and $Y_n$ on the positive end of $M_n^c$.  

\item[{\bf III.}] By Theorem~\ref{theorem:finite} the manifolds $M_n^c$
are geometrically finite hyperbolic cone-manifolds (with no cusps,
since $M$ has no cusps).

\item[{\bf IV.}]  Applying Theorem~\ref{theorem:bilip}, we may deform
the cone-singularities at $\gamma_n^*$ and $\beta_n^*$ back to $2\pi$
fixing the conformal boundary of $M_n^c$ to obtain quasi-Fuchsian
hyperbolic 3-manifolds $Q(X_n,Y_n)$.

\item[{\bf V.}] Since the lengths of $\gamma_n^*$ and $\beta_n^*$ are
tending to zero, the metric distortion of the cone-deformation outside
of tubular neighborhoods of the cone-singularities is tending to zero.
Since the geodesics $\gamma_n^*$ and $\beta_n^*$ are exiting the ends
of $M$, larger and larger compact subsets of $M$ are more and more
nearly isometric to large compact subsets of $Q(X_n,Y_n)$ for $n$
sufficiently large.  Convergence of $Q(X_n,Y_n)$ to $M$ follows.
\end{enumerate}

Next, we detail our approach to the general
doubly-degenerate case, which handles ends with bounded geometry and
ends with arbitrarily short geodesics transparently.

\section{Drilling and asymptotic isolation of ends}
\label{section:drilling}

It is peculiar that manifolds $M \in AH(S)$ of {\em mixed type},
namely, doubly-degenerate manifolds with one bounded
geometry end and one end with arbitrarily short geodesics, present some
recalcitrant difficulties that require new techniques.  Here is an
example of the type of phenomenon that is worrisome:

\bld{Example.}   Consider a sequence $Q(X_n, Y)$ tending to a limit
$Q_\infty$ in the Bers slice $B_Y$ for which $Q_\infty$ is partially
degenerate, and for which $Q_\infty$ has arbitrarily short geodesics.
Allowing $Y$ to vary in Teichm\"uller space, we obtain a {\em limit Bers
slice $B_\infty$ associated to the sequence $\{X_n\}$} (this
terminology was introduced by McMullen \cite{McMullen:graft}).
The limit Bers slice $B_\infty$ is
an embedded copy of $\Teich(S)$ in $AH(S)$ consisting of 
manifolds $$M(Y') = \lim_{n \to \infty} Q(X_n,Y') \ \ \ \text{where} \
\ Y'\ \ \text{lies in} \ \ \Teich(S).$$ Each $M(Y')$ has a degenerate
end that is bi-Lipschitz diffeomorphic to $Q_\infty$ (see, e.g.,
\cite[Prop. 3.1]{McMullen:book:RTM}), but the bi-Lipschitz constant
depends on $Y'$.  

If, for example, $\delta$ is a simple closed curve on $S$ and $\tau^n
 (Y) = Y_n$ is a divergent sequence in $\Teich(S)$ obtained via an
iterated Dehn twist $\tau$ about $\delta$, a subsequence of
$\{M(Y_n)\}_{n = 1}^\infty$ converges to a limit $M_\infty$, but there
is no {\em a priori} reason for the degenerate end of $M_\infty$ to be
bi-Lipschitz diffeomorphic to that of $M(Y)$.  The limiting geometry
of the ends compactified by $Y_n$ could, in principle, bleed over into
the degenerate end, causing its asymptotic structure to change in the
limit.  (We note that such phenomena would violate Thurston's {\em
ending lamination conjecture} since $M_\infty$ has the same ending
lamination associated to its degenerate end as does $M(Y)$).

\bld{Isolation of ends.}  For a convergent sequence of quasi-Fuchsian
manifolds $Q(X_n,Y_n) \to N$, we seek some way to isolate the limiting
geometry of the ends of $Q(X_n,Y_n)$ as $n$ tends to infinity.  Our
strategy is to employ the drilling theorem in a suitably chosen family
of convergent approximates $Q(X_n,Y_n) \to N$ for which a curve
$\gamma$ is short in $Q_n = Q(X_n,Y_n)$ for all $n$.  We prove that
drilling $\gamma$ out of each $Q_n$ to obtain a drilled manifold
$Q_n(\gamma)$ produces a sequence converging to a drilled limit $N(\gamma)$
whose higher genus ends are bi-Lipschitz diffeomorphic to those of $N$.

An application of the covering
theorem of Thurston and Canary
\cite{Thurston:book:GTTM,Canary:inj:radius} then demonstrates that the
limiting geometry of the negative end of $N$ depends only on the
sequence $\{X_n\}$ and the limiting geometry of the positive end of
$N$ depends only on the sequence $\{Y_n\}$.

When $N$ has no such short geodesic $\gamma$, the ends depend only on
the end invariant $\nu(N)$, since in this case $N$ has bounded
geometry and Theorem~\ref{theorem:bounded} applies.  These arguments
are summarized in the following isolation theorem for the asymptotic
geometry of $N$ (cf. Theorem~\ref{theorem:isolation:intro}).

\begin{theorem}{\sc Asymptotic Isolation of Ends}
Let $Q(X_n,Y_n) \in AH(S)$ be a sequence of quasi-Fuchsian manifolds
converging algebraically to the cusp-free limit manifold $N$.  Then,
up to marking and orientation preserving bi-Lipschitz diffeomorphism,
the positive end of $N$ depends only on the sequence $\{ Y_n \}$ and
the negative end of $N$ depends only on the sequence $\{ X_n \}$.
\label{theorem:isolation}
\end{theorem}

We now argue that as a consequence of Theorem~\ref{theorem:isolation}
we need only show that each end of a doubly-degenerate manifold $M$
arises as the end of a singly-degenerate manifold lying in a Bers
boundary.
\begin{defn}
Let $E$ be an  end of a complete hyperbolic
3-manifold $M$.  If $E$ admits a marking and orientation preserving
bi-Lipschitz diffeomorphism to an end $E'$ of a manifold $Q$ lying in
a Bers compactification, we say $E$ is {\em realized in a Bers
compactification by $Q$}.
\label{definition:ends:realizable}
\end{defn}

If, for example, the positive end $E^+$ of $M$ is realized by
$Q_\infty^+$ on the Bers boundary $\bdry B_X$ then there are by
definition surfaces $\{Y_n\}$ so that $Q(X,Y_n)$ converges to
$Q_\infty^+$, so $E^+$ depends only on $\{Y_n\}$ up to bi-Lipschitz
diffeomorphism.  Arguing similarly, if $E^-$ is realized by
$Q_\infty^-$ on the Bers boundary $\bdry{B_Y}$, the approximating
surfaces $\{X_n\}$ for which $Q(X_n,Y)
\to Q_\infty^-$ determine $E^-$ up to bi-Lipschitz diffeomorphism.

By an application of Theorem~\ref{theorem:isolation}, if the manifolds
$Q(X_n,Y_n)$ converge to a cusp-free limit
$N$, then the negative end $E^-_N$ is bi-Lipschitz diffeomorphic to
$E^-$ and the positive end $E^+_N$ is bi-Lipschitz diffeomorphic to $E^+$.
We may glue bi-Lipschitz diffeomorphisms $$\psi^- \colon E^-_N \to E^-
\ \ \
\text{and} \ \ \
\psi^+ \colon E^+_N \to E^+$$
along the remaining compact part to obtain a global bi-Lipschitz
diffeomorphism $$\psi \colon N \to M$$ that is marking and orientation
preserving. Applying Sullivan's rigidity theorem
\cite{Sullivan:linefield}, $\psi$ is homotopic to an isometry, so $Q(X_n,Y_n)$ converges to $M$.

\bld{Realizing ends in Bers compactifications.}  To complete the proof of
Theorem~\ref{theorem:main}, then, we seek to 
realize each end of the doubly-degenerate manifold $M$ on a Bers
boundary; we restate Theorem~\ref{theorem:ends:realizable} here.
\begin{theorem}{\sc Ends are Realizable}
Let $M \in AH(S)$ have no cusps.  Then each end of $M$ is
realized in a Bers compactification.
\end{theorem}

In the case that $M$ has a conformal boundary component $Y$, the theorem
asserts that $M$ lies within the Bers compactification
$\closure{B_Y}$.  This is the main result of \cite{Bromberg:bers},
which demonstrates all such manifolds are limits of quasi-Fuchsian
manifolds.

We are left to attend to the case when $M$ is doubly-degenerate.  As
one might expect, the discussion breaks into cases depending on
whether an end $E$ has bounded geometry or arbitrarily short
geodesics.  We discuss the positive end of $M$; one argues
symmetrically for the negative end.
\begin{enumerate}
\item If a bounded geometry end $E$ has ending lamination $\nu$,
choose a measured lamination $\mu$ with support $\nu$ and a
sequence of weighted simple closed curves $t_n \gamma_n \to \mu$.
Choose $Y_n$ so that $\ell_{Y_n}(\gamma_n) < 1$.

\item If $\gamma_n^*$ are arbitrarily short geodesics exiting 
the end $E$, we apply the drilling theorem to 
${\rm Gr}^\pm(\gamma_0,\gamma_n,M)$ to
send the cone-angles at $\gamma_0^*$ and $\gamma_n^*$ to $2 \pi$.
The result is a sequence $Q(X,Y_n)$ of quasi-Fuchsian manifolds.
\end{enumerate}
We wish to show that after passing to a subsequence $Q(X,Y_n)$ converges to
a limit $Q_\infty$ that realizes $E$ on the Bers boundary $\bdry B_X$.

\ital{Bounded geometry.}  When $E$ has bounded geometry, we employ
\cite{Minsky:KGCC} to argue that its end invariant $\nu$ has {\em
bounded type}.  This condition ensures that any end with $\nu$ as its
end invariant has bounded geometry.  The condition
$\ell_{Y_n}(\gamma_n) < 1$ guarantees that $\ell_{Y_n}(t_n \gamma_n)
\to 0$ so that any limit $Q_\infty$ of $Q(X,Y_n)$ has $\nu$ as its
end-invariant (by \cite{Brock:length}, applying
\cite[Thm. 3]{Bers:bdry}).  We may therefore apply a relative version of
Minsky's ending lamination theorem for bounded geometry (see
\cite{Minsky:ends}, and an extension due to Mosher
\cite{Mosher:elc} that treats the case when the manifold
may not possess a {\em global} lower bound to its injectivity radius)
to conclude that $Q_\infty$ realizes $E$.

\ital{Arbitrarily short geodesics.}
If $E$ has an exiting sequence $\{\gamma_n\}$ of arbitrarily short
geodesics, we argue using Theorem~\ref{theorem:bilip} 
that $Q(X,Y_n)$ converges in the Bers boundary $\bdry B_X$ to a limit
$Q_\infty$ that
realizes $E$.  

\bld{Binding realizations.} As a final detail we mention that 
to apply Theorem~\ref{theorem:isolation}, we require a {\em convergent} sequence
$Q(X_n,Y_n) \to N$ so that the limit $Q^- = \lim Q(X_n,Y_0)$ realizes
the negative end $E^-$ of $M$ and the limit $Q^+ = \lim Q(X_0,Y_n)$
realizes the positive end $E^+$.  

By an application of \cite{Brock:Boundaries}, the realizations
described in our discussion of 
Theorem~\ref{theorem:ends:realizable} produce surfaces $\{X_n\}$ and
$\{Y_n\}$ that converge up to subsequence 
to laminations in Thurston's compactification of Teichm\"uller space
that {\em bind the
surface} $S$.
Thus, an application of  Thurston's {\em double limit
theorem} (see \cite[Thm. 4.1]{Thurston:hype2},
\cite{Otal:book:fibered})
implies that $Q(X_n,Y_n)$ converges to a cusp-free limit $N$ after
passing to a subsequence.

\section{Incompressible ends}
\label{section:general}
We conclude the paper with a brief discussion of the proof of
Theorem~\ref{theorem:main} when $M$ is not homotopy equivalent to a
closed surface.  

Since $M$ has incompressible ends, Theorem~\ref{theorem:tameness}
implies that $M$ is homeomorphic to the interior of a compact
3-manifold $N$.  Equipped with a homotopy equivalence or {\em marking}
$ f \colon N \to M,$ the manifold $M$ determines an element of the
{\em algebraic deformation space} $AH(N)$ consisting of all such
marked hyperbolic 3-manifolds up to isometries preserving orientation
and marking, equipped with the topology of algebraic convergence.

By analogy with the quasi-Fuchsian locus, the subset $AH(N)$
consisting of $M'$ that are geometrically finite, cusp-free and
homeomorphic to $M$ is parameterized by the product of Teichm\"uller
spaces $$
\Teich(\bdry N) =
\prod_{X \subset \bdry N} \Teich(X).
$$ 

In this situation, the cover $\widetilde{M}$ corresponding to an end
$E$ of $M$ lies in $AH(S)$.  Theorem~\ref{theorem:ends:realizable}
guarantees that if $E$ is degenerate it is realized on a Bers
boundary; indeed, since $M$ is cusp-free and $M$ is not homotopy equivalent
to a surface, it follows that $\widetilde{M}$ is itself
singly-degenerate, so Theorem~\ref{theorem:ends:realizable} guarantees
that $\wt M$ lies in a Bers compactification.

The remaining part of Theorem~\ref{theorem:main}, then, follows
from the following version of
Theorem~\ref{theorem:isolation:intro}.
\begin{theorem}{\sc Asymptotic Isolation of Ends II}
Let $M$ be a cusp-free complete hyperbolic 3-manifold with
incompressible ends homeomorphic to $\interior (N)$.  Let $M_n \to M$
in $AH(N)$ be a sequence of cusp-free geometrically finite hyperbolic
manifolds so that each $M_n$ is homeomorphic to $M$.  Let $(E^1,
\ldots, E^m)$ denote the ends of $M$, and let $\bdry M_n = X^1_n
\disjunion\ldots\disjunion X^m_n$ be the corresponding points in
$\Teich(\bdry N)$.  Then, up to 
marking preserving bi-Lipschitz diffeomorphism, $E^j$ depends only on
the sequence $\{X^j_n\}$.
\label{theorem:isolation:II}
\end{theorem}

In the case not already covered by Theorem~\ref{theorem:isolation},
the covers of $M_n$ corresponding to a fixed boundary component are
quasi-Fuchsian manifolds $Q(Y_n,X^j_n)$.  Their limit is the
singly-degenerate cover of $M$ corresponding to $E^j$, so the surfaces
$Y_n$ range in a compact subset of Teichm\"uller space.  

Again, it
follows that the marked bi-Lipschitz diffeomorphism type of the end
$E$ does not depend on the surfaces $Y_n$.  Theorem~\ref{theorem:main}
then follows in this case from an application of
Theorem~\ref{theorem:ends:realizable} to each end degenerate end $E$
of $M$, after an application of Sullivan's rigidity theorem
\cite{Sullivan:linefield}.

%\bibliographystyle{math}
%\bibliography{math}

\noindent{\tiny \sc Department of Mathematics, University of Chicago,
5734 S. University Ave., Chicago, IL 60637}

\smallskip

\noindent{\tiny \sc Department of Mathematics,
California Institute of Technology,
Mailcode 253-37,
Pasadena, CA 91125}

\end{document}